# Mental Geometry

Nicolas Bouleau

**Summary**. This article illustrates pedagogy through training in the handling of abstractions. Mental arithmetic is not limited to numerical calculation; one can mentally calculate primitives and simplify analytical expressions. Even if there is software that does this very well, this training retains its pedagogical value. Can we go further and consider geometric mental arithmetic: mentally proceeding with transformations of simple figures allowing the calculation of areas or volumes?

It turns out that the intuition that allowed Archimedes to obtain his main geometric results, if we take only the ideas without the old-fashioned style, provides the opportunity for a pleasant and quite rich mental game that I present here in the form of a short narrative dialogue, not a philosophical tale because it does not bring any thesis, simply a story to be classified among the invitations to exercise the mind. It starts with the area of a triangle and ends with Guldin's two theorems.

To show that you don't need a blackboard for this kind of exercise, let's imagine that the setting is a group of pupils on a long-distance hiking trail with their maths, history and philosophy teachers.

*The maths teacher* – Why doesn't the area of a triangle change when you move its vertex parallel to the opposite side?

*A student* – Because the area is half the product of the side and the height.

*The maths teacher* – Here is an explanation indeed, but is it a true understanding? We can understand the invariance of the area during the movement. What do you think of the following argument: let's cut the triangle into thin slices parallel to the fixed side, when the vertex moves, the slices slide over each other. We can see that the area is preserved.

*Another pupil* – But there are little steps at the edges that cause a problem.

*The maths teacher* – Yes, however it is easy to convince oneself that these problems are less serious the thinner the slices and that, in the limit, the reasoning is perfectly rigorous. It is this intuition that I would like to familiarise you with.

*The history teacher* – Isn't it because of this pushing to the limit that Archimedes uses his method of exhaustion, the famous method mentioned many times by the Latins, the Arabs and during the Renaissance?

*The maths teacher* – Indeed, Archimedes reasons as follows: suppose that the area of the deformed triangle exceeds the area of the original triangle by epsilon. By cutting the latter into slices thin enough that the sum of the areas of the external rectangles minus that of the areas of the internal rectangles is less than epsilon, then by moving these rectangular slices to obtain a frame for the deformed triangle, we obtain a contradiction in the inequalities.
It is therefore the method of exhaustion that ensures our passage to the limit. We can always refer to it in case of doubt, but the most interesting thing is to develop our intuition because it is this that enriches our capacity to act. So I ask a new question. Why is the area of a circle equal to half the circumference multiplied by the radius?

*A pupil* – I guess the answer is not that it is $\pi r^2$.

*Another pupil* – Because if we cut the circle like a pie into thin slices, we can unfold it by letting the crusts of the neighboring slices touch and align along a line segment equal to the circumference. The segments then form a sawtooth pattern. If we imagine bringing all the points back to one of them parallelly to the line segment, as we have just seen with triangles, the area is unchanged, we are left with a single triangle whose area is the area we are looking for.

*The maths teacher* – Very good, and you feel that the mistake made disappears if the parts are increasingly thin.

The group set off again. With a certain impatience to hear what came next, we chatted. The history and philosophy teachers mentioned Euclid, Eratosthenes and the Alexandrian School, the birth of Greek science and the irrationals. At the next stop, the maths teacher resumed his talk.

*The maths teacher* – Now let's talk about volumes. Who can tell me why a twisted column, like the ones you see in baroque churches, has the same volume as a straight column with the same base and height?

*A pupil* – It's like the piles of chocolate bars in sweet shops.

*The maths teacher* – Can you explain?

*The student* - If the horizontal sections of the truncated column are circles identical at the base, it is obvious from the reasoning of the slices, just as from a stack of coins one can obtain a helical stack of the same volume and height by shifting them.

*The maths teacher* – And the volume of a cone? It is now clear that it does not change if we move its vertex in a plane parallel to the base, and that moreover it does not depend on the shape of the base but on its area, as we can see by cutting it into small parts and considering the corresponding small cones. Therefore, that this volume is one third of the area of the base times the height results from any particular case, for example from the fact that a cube can be divided from its centre into six equal pyramids.
     We also see that by dividing a sphere into small cones emanating from its centre and opening it up like a mango or a papaya and then bringing all the points together into a single one, the volume of the sphere is equal to one third of its surface area by its radius. Similarly, a regular polyhedron circumscribed about a sphere is composed of cones emanating from the centre of the sphere, and its volume is therefore equal to one third of its surface area multiplied by the radius of the sphere. It is interesting to note that this property is still true if the polyhedron is irregular, provided that all its faces are tangent to the sphere, even if the number of faces is infinite, as with the cylinder circumscribed about the sphere, for example.

*The history teacher* – You are very close to the results of the treatise on the sphere and cylinder of Archimedes, of which he was most proud, if Plutarch is to be believed.

*The maths teacher* – Our intuition can now easily obtain them. Let us consider the cylinder circumscribed around a sphere and cut it with slices parallel to the bases. It turns out that in each slice the surface area of the cut-out sphere is equal to the area of the cut-out portion of the cylinder. This is immediate with a little calculation if the slice is infinitely thin. So it is true, by exhaustion, if the slice is arbitrary. Hence it follows that the area of the sphere is equal to the lateral area of the cylinder, which is two-thirds of the total area of the cylinder, including the bases; from this, and from what has been said about the volume of the sphere and that of circumscribed polyhedron, it follows that the volume of the sphere is two-thirds of the volume of the cylinder.

*The philosophy teacher* – All this does not prove thought without words, because we are certainly following you, but you do use many words to lead us through your constructions.

*The history teacher* - You are reviving the so-called 'indivisibles' method of the seventeenth century, named by the Italian Cavalieri and perfected by Roberval, but it seems to me that with Fermat, Descartes, Leibniz and then Newton, these methods were supplanted by those of differential calculus. What is the point of revisiting outdated ideas?

*The maths teacher* – I'm showing you that you can do mental arithmetic without numbers, using geometry alone. Integral calculus has killed the intuition of indivisibles, it's a shame, there was something interesting there. As for proving thought without words, I can see no other way now than to let the intuition I have aroused among the pupils run free to see if this thought develops by itself.

The group set off again, but this time in silence. Nobody dared to speak. Already the light was fading in the east and the first shadows were creeping in, taking on unexpected shapes. After an hour, raised voices could be heard here and there.

*A group of pupils:* Sir, we have the volume of a portion of a cylinder cut by two planes.

*The maths teacher* – We're listening.

*One of the pupils* – Here we go. Let's consider a sphere and divide it up with meridional planes like an orange or time zones if you prefer. On each zone, the equator draws an arc. Let's spread out our orange on a table so that the ends of the equator touch the table and line up. Then let's mentally deform these quarters by moving the segments that make up their ribs, translating them parallel to the plane of the table to bring them back to one of them. The lateral surface is preserved as well as the total volume. We obtain a portion of a right cylinder cut by two planes passing through a perpendicular to its axis, I don't know if it has a name.

*The maths teacher* – Good. So we see that if we cut a vertical cylinder of revolution with a plane passing through a diameter of the base, the sector, which is shaped like a mitre, has its lateral surface proportional to its total height and the same goes for its volume, and your reasoning provides the two coefficients of proportionality.

*The philosophy teacher* – Excuse me, I'd like to see a diagram. Geometry has always been done with diagrams. It is normal to express oneself with drawings, they are the ancestral forms of the symbols of language. Architects and engineers make drawings, painters sometimes do too. Can't we say that figures are to geometry what formulas are to arithmetic or algebra?

*The maths teacher* – Not quite. Formulas are in the realm of the discrete, figures in that of the continuous, there are questions of precision that can be misleading, apparent implications that are not mathematical implications. Geometry is the art of deductive thinking about lines. When you manage to handle them, I believe that mental figures are the only ones that are logically perfect.

*The student* - My reasoning is based solely on the deformation of an orange wedge. Let's assume that such a wedge is a section of a cylinder by two planes. We simply use the fact that neither the volume nor the lateral surface change if the edge is moved parallel to the axis of the cylinder. This is a generalisation of the property of the area of a triangle not to change when a vertex is moved parallel to the opposite side.

The philosophy teacher said nothing. He noted that the student's intuition had produced knowledge of the same nature, it seemed, as science, awkwardly expressed by the words of ordinary language, but whose rigour was in no way inferior to the most accomplished texts.

This would mean that thinking could be a spontaneous, uncontrolled activity. He asked the students what had sparked their idea. Because of the heat, one of them replied, we've been getting more and more urges to eat oranges...

*The history teacher* – This section of a cylinder cut by two planes passing through a perpendicular to its axis is sometimes called Archimedes' clog, he was quite proud to have found its volume and uses this calculation as an example to demonstrate his method to Eratosthenes, it was the first 'kybismos', the exact measurement of a volume bounded in part by a curved surface. But your approach is simpler than his, which uses two successive equilibriums, one of which results from very ingenious metric relationships. He also confirms his result by another method based on the area of the section of the parabola, which he had obtained previously. Finally, you made me grasp his thinking much better than the comments of Eutocius, which, I must say, I never understood at all.

*The maths teacher* – Eutocius is a hardworking man who, seven centuries after Archimedes, repeats his reasoning, weighing down everything that was conceptual and living. It's uninteresting. The same remark is made by Jean Dieudonné about the comments on the work of Evariste Galois by some of his contemporaries. They are complicated and convoluted where Galois' text seems simple and clear to us. Archimedes himself said how he wanted to be read:

> I am convinced, in fact, that researchers, either of our time or of the future, will find, by applying the method that I will have made known, still other propositions that did not come to my mind.

I must point out, however, that today we have only glimpsed part of Archimedes' heuristic thinking, which consists of deforming infinitesimal elements while preserving the areas in the plane or the volumes in space, but Archimedes complements this tool with another one, which is to deform them while preserving the balance in relation to an axis of the supposedly heavy parts. This is also very fruitful, of course.

To address slightly more scholarly questions, but in the same vein, it is reasonable to assume that the students returned to a classroom and sat down to catch their breath.

*The maths teacher* – If we use the Archimedes method to its full potential, with the equilibrium reasoning that I mentioned earlier, we can easily obtain Guldin's two theorems relating to bodies of revolution, which Pappus was familiar with from the fourth century.[1] First, it is sufficient to demonstrate the following lemma: let there be a closed curve in the horizontal plane serving as the base of a straight vertical cylinder infinite at the top and bottom. If we cut this cylinder with a plane passing through the centre of gravity of the base, the two portions of the cylinder cut above the horizontal plane (between the base and the cutting plane) and below the horizontal plane (between the cutting plane and the base) have the same volume. And similarly, if the cylinder is cut by a plane passing through the centre of gravity of the boundary curve of this base, the lateral surfaces of the cylinder portions above and below the horizontal plane are equal.

     Let us consider the volumes. Note that, with respect to a line in its plane passing through the centre of gravity, the supposedly weighty base-surface is in equilibrium. This equilibrium is an equality of moments for the surface elements of the base, and the numerical value of the moment of a surface element is, to within a coefficient, the volume of the portion above this element or below it, as the case may be, limited by the oblique plane.

---

[1] Paul Guldin (1577-1643) is traditionally credited with the discovery of these theorems (which we state below), but there is a text by Pappus of Alexandria in which the theorem relating to volumes is stated. Pappus claims to have a proof of it, which may have appeared in the volumes of his work that have not survived.

The same reasoning applies to lateral surfaces. The lemma is simpler than many of Archimedes' arguments. The stated equalities directly translate into equalities between moments without having to make any distortions.

Then we consider a body of revolution whose section by a meridian half-plane has any shape. Using the same method of slices by planes passing through the axis, unfolding, and then regrouping, the body is transformed into a right cylinder with the section as its base, limited at the other end by an oblique plane and with the same volume and lateral area as the initial body. If the base is horizontal and the cylinder vertical, each vertical fibre of the cylinder has the length of the circumference described by the point where it meets the base. By applying the previous lemma, the volume of the cylinder is not changed by cutting it horizontally at the height of the circumference described by the centre of gravity of the section. We have Guldin's first theorem: the volume is the product of the area of the section by the circumference described by the centre of gravity. And the same reasoning gives us I second theorem: the lateral area is the product of the length of the boundary curve of the section by the circumference described by the centre of gravity of this curve.

It can be noted that the calculation of indivisibles, that is to say what we have just practised, was in full swing in the seventeenth century with Roberval and the area of the cycloid arch when it was supplanted by integral calculus to the point of being abandoned. However, this geometric method of transforming infinitesimal elements to prove the equality of two areas could be applied *a priori*, even if the integrals that provide these areas are incalculable. The two methods therefore did not have the same field of application.
The discussions among the pupils were lively. The philosophy teacher took the references from Hadamard's book and, as a final word, proposed a quotation from the philosopher Francis Bacon from the very beginning of the seventeenth century:

> Pure mathematics cure many intellectual defects. If the wit be loo dull, they sharpen it; if too wandering, they fix it; if too inherent in the sense, they abstract it. So that as tennis is a game of no use in itself, but of great use in respect it maketh a quick eye and a body ready to put itself into all postures; so in the mathematics, that use which is collateral and intervenient is no less worthy than that which, is principal and intended. [*Analysis of the advancement of learning*]